\documentclass[a4,11pt]{article}
\usepackage{amsfonts}
\usepackage{amsmath}
\usepackage{amssymb}
\usepackage{amsthm}
\theoremstyle{plain}
 \newtheorem{theorem}{Theorem}[section]
 
 \newtheorem*{theorem*}{Theorem}
 
 \newtheorem{lemma}[theorem]{Lemma}
 \newtheorem{corollary}[theorem]{Corollary}

\theoremstyle{definition}
 \newtheorem{definition}[theorem]{Definition}
\theoremstyle{remark}
 \newtheorem{remark}[theorem]{Remark}
 
\numberwithin{equation}{section}

\title{Characterizations of regularity for certain $Q$-polynomial association schemes}

\author{Sho Suda\\
{\small Division of Mathematics, Graduated School of Information Sciences, Tohoku University,} \\
{\small 6-3-09 Aramaki-Aza-Aoba, Aoba-ku, Sendai 980-8579, Japan}\\
{\small suda@ims.is.tohoku.ac.jp}}
\date{\today}
\begin{document}
\maketitle

\begin{abstract}
In previous paper \cite{S}, it is shown that linked systems of symmetric designs with $a_1^*=0$ and mutually unbiased bases (MUB) are triply regular association schemes.
In this paper, we characterize triple regularity of linked systems of symmetric designs by its Krein number.
And we prove that maximal MUB carries a quadruply regular association scheme and characterize the quadruple regularity of MUB by its parameter.     
\end{abstract}
\section{Introduction}
We study the regularity of $3$-class (respectively $4$-class) $Q$-polynomial association schemes with $Q$-antipodal (respectively both $Q$-antipodal and $Q$-bipartite).

In Section~\ref{LSD}, we consider linked systems of symmetric designs.
Systems of projective designs, that were defined by P. J. Cameron \cite{C}, 
are the combinatorial object of 
finite doubly transitive groups which have more than two pairwise inequivalent permutation  representations with the same permutation character. 
We call it linked systems of symmetric designs,
if symmetric designs appearing in systems of projective designs are all same  parameters. 
Noda \cite{N} showed several inequalities concerning the parameters of linked systems of  symmetric designs.  
Mathon \cite{M} showed every linked system of symmetric designs carries a $3$-class association 
scheme and calculated its eigenmatrices. 
It implies that these association schemes are $Q$-polynomial with $Q$-antipodal.
Conversely van Dam \cite{D} showed every $3$-class $Q$-polynomial association scheme 
with $Q$-antipodal arises from a linked system of symmetric designs. 
The author \cite{S} proved that 
every linked system of symmetric designs with $a_1^*=0$ is a triply regular association scheme.
Main theorem in this section is the converse proposition, that is, 
if a linked system of symmetric designs is triply regular, then $a_1^*=0$.   
This proof is essentially due to \cite[Theorem 2]{N}.

In Section~\ref{qua}, we consider the quadruple regularity of symmetric association schemes.
We define the quadruple regularity and 
give the sufficient condition that spherical designs become the quadruply regular symmetric association schemes. 

In Section~\ref{MUB}, we consider the real mutually unbiased bases (MUB).
One important problem of real MUB is to determine the maximal number of real MUB in $\mathbb{R}^d$.
It is well known that its number is at most $d/2+1$. 
Real MUB is said to be maximal if equality holds.
Recently W. J. Martin et al. \cite{LMO} showed that there is a one-to-one correspondence between 
real MUB and $4$-class $Q$-polynomial association schemes which is both $Q$-bipartite and 
$Q$-antipodal.  
Moreover W.J. Martin et al. \cite{MMW} had shown that 
a $4$-class $Q$-polynomial association schemes which is both $Q$-bipartite and 
$Q$-antipodal is obtained by the extended $Q$-bipartite double of
a linked system of symmetric design with certain parameters.  
The author proved in \cite{S} that every MUB carries a triply regular association scheme.
The main theorem in this section is that MUB carries a quadruply regular association scheme
if and only if MUB is maximal.
\section{Linked systems of symmetric designs}\label{LSD}
\begin{definition}
Let $(X_i,X_j,I_{i,j})$ be an incidence structure satisfying 
$X_i\cap X_j=\emptyset$, $I_{j,i}^t=I_{i,j}$ for any distinct integers $i,j \in \{1,\dots,f\}$. 
We put $X=\bigcup_{i=1}^f X_i$,
$I=\bigcup_{i\neq j} I_{i,j}$.
$(X, I)$ is called a linked system of symmetric $(v,k,\lambda)$ designs if the following conditions hold:
\begin{enumerate}
\item for any distinct integers $i,j\in\{1,\dots,f\}$, $(X_i,X_j,I_{i,j})$ is a symmetric $(v,k,\lambda)$ design,
\item for any distinct integers $i,j,l \in\{1,\dots,f\}$, and for any $x\in X_i, y\in X_j$, the number of $z\in X_l$ incident with 
both $x$ and $y$ depends only on whether $x$ and $y$ are incident or not, and does not depend on $i,j,l$.
\end{enumerate}
\end{definition}
We define the integers $\sigma,\tau$ by
\[
|\{z\in X_l \mid (x,z)\in I_{i,l}, (y,z)\in I_{j,l}\}|=
\begin{cases}\sigma&\text{ if }\ (x,y)\in I_{i,j},\\
\tau&\text{ if }\ (x,y)\not\in I_{i,j},
\end{cases}
\]
where $i,j,l \in\{1,\dots,f\}$ are distinct and $x\in X_i$, $y\in X_j$.
Theorem 1 in \cite{C} shows 
$$(\sigma,\tau)=\left(\frac{1}{v}(k^2\mp\sqrt{n}(v-k)),\frac{k}{v}(k\pm\sqrt{n})\right),$$
where $n=k-\lambda$.
Considering complement designs $(X_i,X_j,\overline{I_{i,j}})$ for any distinct integers $i,j\in \{1,\ldots,f\}$, 
we can assume either $(\sigma,\tau)=(\frac{1}{v}(k^2-\sqrt{n}(v-k)),\frac{k}{v}(k+\sqrt{n}))$ or $(\frac{1}{v}(k^2+\sqrt{n}(v-k)),\frac{k}{v}(k-\sqrt{n}))$.

We obtain a $Q$-antipodal $3$-class $Q$-polynomial association scheme $(X,\{R_i\}_{i=0}^3)$ where
\begin{align*}
R_0&= \{(x,x) \mid x\in X\} ,\\
R_1&= \{(x,y)\mid x\in X_i,y\in X_j, (x,y)\in I_{i,j}\ \text{ for some }\ i\neq j\} ,\\
R_2&= \{(x,y)\mid x,y\in X_i, x\neq y \text{ for some }\ i\} ,\\
R_3&= \{(x,y)\mid x\in X_i,y\in X_j, (x,y)\not\in I_{i,j}\ \text{ for some }\ i\neq j\}.
\end{align*}
Conversely every $Q$-antipodal $3$-class $Q$-polynomial association scheme with equivalence relation $R_0\cup R_2$ arises 
from a linked system of symmetric designs in \cite[Theorem 5.8]{D}.

\begin{theorem}\label{triLSD}
$(X,\{R_i\}_{i=0}^3) $ is a $Q$-polynomial association scheme which is $Q$-antipodal  with equivalence relation $R_0\cup R_2$.
Then the following are equivalent.
\begin{enumerate}
\item $(X,\{R_i\}_{i=0}^3)$ is triply regular. 
\item $a_1^*=0$.
\end{enumerate}
\end{theorem}
\begin{proof}
(2)$\Rightarrow$(1): Follows from \cite[Corollary 6.2]{S}.

(1)$\Rightarrow$(2):Let $\{X_i,\ldots,X_f\}$ be a system of imprimitivity with respect to the equivalence relation $R_0\cup R_2$ and $(X,R_1)$ a linked system of symmetric $(v,k,\lambda)$ designs. 
Assume that 
$$ \sigma=\frac{1}{v}(k^2-\sqrt{n}(v-k)),\quad \tau=\frac{k}{v}(k+\sqrt{n}). $$
By the assumption of triple regularity, the following number 
$$|\{w\in X \mid (x,w),(y,w),(z,w)\in R_1\}|$$
for distinct points $x,y,z\in X_1$ does not depend on $x,y,z\in X_1$. 
This implies that a pair $(X_1,\bigcup_{i=2}^fX_i)$ is a $3$-design, therefore equality holds in \cite[Theorem 2]{N}.
It follows that
$$f-1=\frac{(v-2)\sqrt{k(v-k)}}{(v-2k)\sqrt{v-1}}.$$ 
This implies $a_1^*=0$ (See \cite[p.14]{S}).
\end{proof}
\begin{remark}
Mathon \cite{M} pointed out that the inequality in \cite[Theorem 2]{N} is equivalent to $a_1^*\geq0$.
\end{remark}
\section{Quadruple regularity of symmetric association schemes}\label{qua} 
\begin{definition}
Let $(X,\{R_i\}_{i=0}^d)$ be a symmetric association scheme. Then the association scheme $X$ is said to be quadruply regular if,
for all $I=(i_1,i_2,i_3,i_4)\subset \{0,1,\ldots,d\}^4$, $J=(j_{\alpha,\beta})_{1\leq \alpha<\beta\leq 4} \subset \{0,1,\ldots,d\}^6$ and $x_1,\ldots,x_4\in X$ such that $(x_k,x_l)\in R_{j_{k,l}}$ for any $1\leq k<l\leq4$, 
the number
$$ |R_{i_1}(x_1)\cap R_{i_2}(x_2)\cap R_{i_3}(x_3)\cap R_{i_4}(x_4)| $$
depends only on $I$, $J$ and not on $x_1,\ldots,x_4$.
\end{definition}

Let $(X,\{R_i\}_{i=0}^d)$ be a symmetric association scheme.
We define the $i$-th subconstituent with respect to $z\in X$ by $R_i(z):=\{y\in X\mid (z,y)\in R_i\}$ and the $(i,j)$-th subconstituent with respect to $(z_1,z_2)\in X\times X$ by $R_{i,j}(z_1,z_2):=R_i(z_1)\cap R_j(z_2)$.
We denote by $R_{i,j,k,l}^{m,n}(z_1,z_2)$ the restriction $R_n$ to $R_{i,j}(z_1,z_2)\times R_{k,l}(z_1,z_2)$ for $(z_1,z_2)\in R_m$. 
Moreover let $(X,\{R_i\}_{i=0}^d)$ be triply regular. 
We denote $p_{l,m,n}^{i,j,k}=|R_m(x)\cap R_n(y)\cap R_l(z)|$
for $x,y,z \in X$ such that $(x,y)\in R_i, (y,z) \in R_j, (z,x)\in R_k$.
Quadruple regularity is characterized by the concept of coherent configuration. 
We omit easy proof of the following lemma.
\begin{lemma}\label{symqua}
A symmetric association scheme $(X,\{R_i\}_{i=0}^d)$
is quadruply regular if and only if $(X,\{R_i\}_{i=0}^d)$ is triply regular and
 for all $m\in\{1,\ldots,d\}$ and $z_1,z_2\in X$ with $(z_1,z_2)\in R_m$, 
$(\bigcup\nolimits_{i,j=1}^dR_{i,j}(z_1,z_2), \{R_{i,j,k,l}^{m,n}(z_1,z_2)\mid 1\leq i,j,k,l\leq d,p_{i,j,m}^{l,k,n}\neq0\})$ 
is a coherent configuration whose parameters depend only on $m$, not on the choice of $z_1,z_2$ with $(z_1,z_2)\in R_m$. 
\end{lemma}
Let $X$ be a finite subset in $S^{d-1}$ with degree $s$, and $A(X)=\{\alpha_1,\ldots,\alpha_s\}$.
For $z_1,z_2\in X$ with $\langle z_1,z_2\rangle=\alpha_m\neq\pm1$,
$X_{i,j}^m=X_{i,j}^m(z_1,z_2)$ will denote the orthogonal projection of 
$\{y\in X\mid \langle  y,z_1\rangle=\alpha_i ,\langle y,z_2\rangle=\alpha_j\}$ 
to $\langle z_1,z_2\rangle^\perp=\{y\in \mathbb{R}^d \mid \langle y,z_1\rangle =\langle y,z_1\rangle=0 \}$, rescaled to lie in $S^{d-3}$.
If $\langle x,z_1\rangle=\alpha_i$, $\langle x,z_2\rangle=\alpha_j$, $\langle  y,z_1\rangle=\alpha_k$, $\langle y,z_2\rangle=\alpha_l$ and $\langle x,y\rangle=\alpha_n$,
then the inner product of the orthogonal projections of $x,y$ to 
$\langle z_1,z_2\rangle^\perp$ rescaled to lie in $S^{d-3}$ is
$$\alpha_{i,j,k,l}^{m,n}:= \frac{(\alpha_n-\alpha_i\alpha_k)(1-\alpha_m^2)-(\alpha_j-\alpha_i\alpha_m)(\alpha_l-\alpha_k\alpha_m)}{\sqrt{(1-\alpha_i^2-\alpha_j^2-\alpha_m^2+2\alpha_i\alpha_j\alpha_m)(1-\alpha_k^2-\alpha_l^2-\alpha_m^2+2\alpha_k\alpha_l\alpha_m)}}.$$
We denote $p_{\alpha,\beta}^{(i,j,m)}(x,y)=|\{z\in X_{i,j}^m\mid \langle x,z\rangle=\alpha,\langle y,z\rangle=\beta\}|$.
\begin{lemma}
Let $X\subset S^{d-1}$ be a finite set and $A'(X)=\{ \alpha_1,\dots,\alpha_s\}$.
Assume that $(X,\{R_k\}_{k=0}^s)$ is a symmetric association scheme,
where $R_k=\{(x,y)\in X\times X \mid \langle x,y\rangle=\alpha_k\}$
$(0\leq k\leq s)$ and $\alpha_0=1$.
Then $|\{(i,j)\in \{1,\ldots,s\}^2\mid X_{i,j}^m(z_1,z_2)\neq \emptyset\}|=|\{(i,j)\in \{1,\ldots,s\}^2\mid p_{i,j}^m\neq0\}|$ for $\langle z_1,z_2\rangle=\alpha_m$. 
\end{lemma}
\begin{proof}
Immediate from definition.
\end{proof}

The following theorem is used to prove Corollary~\ref{quadruple}.
\begin{theorem}[{\cite[Theorem 2.6]{S}}]\label{coherent}
Let $X_i \subset S^{d-1}$ be a spherical $t_i$-design for $i \in \{1,\dots,n\}$.
Assume that $X_i \cap X_j=\emptyset$ or $X_i=X_j$, and $X_i \cap (-X_j)=\emptyset$ or $X_i=-X_j$ for $i ,j \in \{1,\dots,n\}$.
Let $s_{i,j}=|A(X_i,X_j)|$, $s_{i,j}^*=|A'(X_i,X_j)|$ and $A(X_i,X_j)=\{\alpha_{i,j}^1,\ldots,\alpha_{i,j}^{s_{i,j}}\}$, $\alpha_{i,j}^0=1$, when $-1\in A'(X_i,X_j)$, we define  $\alpha_{i,j}^{s_{i,j}^*}=-1$. 
We define $R_{i,j}^k=\{(x,y)\in X_i \times X_j\mid \langle x,y \rangle =\alpha_{i,j}^k \} $.
If one of the following holds depending on the choice of $i,j,k \in \{1, \dots ,n \}$:
\begin{enumerate}
\item $s_{i,j}+s_{j,k}-2 \leq t_j$,
\item $s_{i,j}+s_{j,k}-3 = t_j$ and for any $\gamma \in A(X_i,X_k)$ there exist $\alpha\in A(X_i,X_j),\beta\in A(X_j,X_k)$ such that the number $p_{\alpha,\beta}^j(x,y)$ is independent of the choice of $x\in X_i,y\in X_k$ with $\gamma=\langle x,y\rangle$,
\item $s_{i,j}+s_{j,k}-4 = t_j$ and for any $\gamma \in A(X_i,X_k)$  there exist $\alpha,\alpha'\in A(X_i,X_j),\beta,\beta'\in A(X_j,X_k)$  such that $\alpha\neq \alpha'$, $\beta\neq \beta'$ and 
the numbers $p_{\alpha,\beta}^j(x,y)$, $p_{\alpha,\beta'}^j(x,y)$ and $p_{\alpha',\beta}^j(x,y)$ are independent of the choice of $x\in X_i,y\in X_k$ with $\gamma=\langle x,y\rangle$,
\end{enumerate}   
then $(\coprod \nolimits_{i=1}^n X_i, \{R_{i,j}^k\mid 1\leq i,j\leq n,1-\delta_{X_i,X_j}\leq k\leq s_{i,j}^*\}) $ is a coherent configuration.
The parameters of this coherent configuration are determined by $A(X_i,X_j)$, $|X_i|$, $t_i$, $\delta_{X_i,X_j}$, $\delta_{X_i,-X_j}$, and when $s_{i,j}+s_{j,k}-3=t_j$ (resp. $s_{i,j}+s_{j,k}-4=t_j$), the numbers $p_{\alpha,\beta}^j(x,y)$ (resp. $p_{\alpha,\beta}^j(x,y)$, $p_{\alpha',\beta}^j(x,y)$, $p_{\alpha,\beta'}^j(x,y)$) which are assumed be independent of $(x,y)$ with $\langle x,y\rangle=\gamma$. 
\end{theorem}
The following lemma shows the antipodal double cover of coeherent configurations are also coherent configurations.   
\begin{lemma}\label{anticc}
Let $X_i^+,X_i^-\subset S^{d-1}$ be a finite subset such that $X_i^+=-X_i^-$ for $i\in\{1,\ldots,n\}$.
If $\{X_i^+\}_{i=1}^n$ carries a coherent configuration, 
then $\{X_i^+,X_i^-\}_{i=1}^n$ carries also a coherent configuration.
\end{lemma}
\begin{proof}
We define $X_i^\varepsilon(x,\alpha)=\{w\in X_i^\varepsilon\mid \langle x,w\rangle=\alpha\}$, 
and $X_i^\varepsilon(x,\alpha;y,\beta)=X_i^\varepsilon(x,\alpha)\cap X_i^\varepsilon(y,\beta)$ for $x\in S^{d-1}$, $\varepsilon=+$ or $-$.
Then the following equalities hold:
\begin{enumerate}
\item $X_i^+(x,-\alpha)=X_i^+(-x,\alpha)$,
\item $X_i^+(x,\alpha)=-X_i^-(x,-\alpha)$.
\end{enumerate}
By (1), $X_i^+(x,\alpha;y,\beta)=X_i^+(-x,-\alpha;y,\beta)=X_i^+(x,\alpha;-y,-\beta)$ holds.
By (2), $X_i^+(x,\alpha;y,\beta)=-X_i^-(x,-\alpha;y,-\beta)$ holds.
Therefore $|X_i^+(x,\alpha;y,\beta)|=|X_i^+(-x,-\alpha;y,\beta)|=|X_i^+(x,\alpha;-y,-\beta)|=|X_i^-(x,-\alpha;y,-\beta)|$ holds.
This implies that intersection numbers on $\{X_i^+,X_i^-\}_{i=1}^n$ is determined by the  coherent configuration $\{X_i^+\}_{i=1}^n$. 
\end{proof}
The following corollary gives the sufficient condition of the quadruple regularity of triply regular association schemes obtained from an antipodal finite subset of sphere.
Its proof follows from the same argument of \cite[Corollary 2.9]{S}.
\begin{corollary}\label{quadruple}
Let $X\subset S^{d-1}$ be an antipodal finite subset and $A'(X)=\{ \alpha_1,\dots,\alpha_s\}$ with $\alpha_1>\dots>\alpha_s=-1$.
Assume that $(X,\{R_k\}_{k=0}^s)$ is a triply regular symmetric association scheme,
where $R_k=\{(x,y)\in X\times X \mid \langle x,y\rangle=\alpha_k\}$
$(0\leq k\leq s)$ and $\alpha_0=1$.
Then for $1\leq i,j,k,l,m\leq s-1$ such that $p_{i,j}^m\neq0$ and $p_{k,l}^m\neq0$
\begin{enumerate}
\item $A(X_{i,j}^m(z_1,z_2),X_{k,l}^m(z_1,z_2))=\{
\alpha_{i,j,k,l}^{m,n}
 \mid 0\leq n\leq s, p_{i,j,m}^{l,k,n}\neq0, \alpha_{i,j,k,l}^{m,n}\neq\pm1\}$.
\item $X_{i,j}^m(z_1,z_2)=X_{k,l}^m(z_1,z_2)$ or $X_{i,j}^m(z_1,z_2)\cap X_{k,l}^m(z_1,z_2)=\emptyset$,
and $X_{i,j}^m(z_1,z_2)=-X_{k,l}^m(z_1,z_2)$ or $X_{i,j}^m(z_1,z_2)\cap -X_{k,l}^m(z_1,z_2)=\emptyset$
for any $z_1,z_2\in X$ with $\langle z_1,z_2\rangle=\alpha_m$
And $\delta_{X_{i,j}^m(z_1,z_2),X_{k,l}^m(z_1,z_2)}$, $\delta_{X_{i,j}^m(z_1,z_2),-X_{k,l}^m(z_1,z_2)}$ are independent of $z_1,z_2\in X$ with $\alpha_m=\langle z_1,z_2\rangle$.
\item $X_{i,j}^m(z_1,z_2)$ has the same strength for all $z_1,z_2\in X$ with $\alpha_m=\langle z_1,z_2\rangle$.
\end{enumerate}
Moreover if the assumption $(1),(2)$ or $(3)$ of Theorem~\ref{coherent}
is satisfied for $\{X_{i,j}^m(z_1,z_2)\mid 1\leq i\leq \frac{s-1}{2},1\leq j\leq s-1, p_{i,j}^m\neq0\} \cup\{X_{i,j}^m(z_1,z_2)\mid \frac{s-1}{2}\leq i\leq \frac{s+1}{2},1\leq j\leq\frac{s+1}{2}, p_{i,j}^m\neq0\}$ with $m\neq0$ or $s$,
and when $((i,j),(k,l),(m,n))$ satisfies $(2)$ (resp. $(3)$) the numbers $p_{\alpha,\beta}^{(k,l,m)}(x,y)$ (resp. $p_{\alpha,\beta}^{(k,l,m)}(x,y)$, $p_{\alpha,\beta'}^{(k,l,m)}(x,y)$, $p_{\alpha',\beta}^{(k,l,m)}(x,y)$) which are assumed to be independent of $(x,y)$ with $\gamma=\langle x,y\rangle$ are independent of the choice of $z_1,z_2$ with $\alpha_m=\langle z_1,z_2\rangle$, 
then $(X,\{R_k\}_{k=0}^s)$ is a quadruply regular association scheme.
\end{corollary}
\begin{proof}
(1), (2), (3) follow from arguments similar to that in \cite[Corollary 2.9]{S}.

Fix $z_1,z_2\in X$ with $\alpha_m=\langle z_1,z_2\rangle$.

If $m=0$ or $s$, then $\bigcup_{i,j=1}^sR_i(z_1)\cap R_j(z_2)=\bigcup_{i=1}^sR_i(z_1)$.
The triple regularity of $(X,\{R_k\}_{k=0}^s)$ is equivalent that $\bigcup_{i,j=1}^sR_i(z_1)\cap R_j(z_2)$ to be a coherent configuration whose parameters are independent of $z_1,z_2$ with $\langle z_1,z_2\rangle=\pm1$.

If $1\leq m\leq s-1$, then $X_{i,s}^m(z_1,z_2)\neq \emptyset$ if and only if $X_{s,i}^m(z_1,z_2)\neq \emptyset$ if and only if $i=s-m$ hold, and then
 $X_{s,m-s}^m(z_1,z_2)=\{-z_1\}$, $X_{s-m,s}^m(z_1,z_2)=\{-z_2\}$ hold. 
Moreover $X_{i,j}^m=-X_{s-i,s-j}^m$ hold for any $1\leq i,j\leq s-1$.
By Lemma~\ref{anticc}, it is sufficient to show that $\{X_{i,j}^m(z_1,z_2)\mid 1\leq i\leq \frac{s-1}{2},1\leq j\leq s-1, p_{i,j}^m\neq0\} \cup\{X_{i,j}^m(z_1,z_2)\mid \frac{s-1}{2}\leq i\leq \frac{s+1}{2},1\leq j\leq\frac{s+1}{2}, p_{i,j}^m\neq0\}$ carries a coherent configuration 
whose parameters are independent of $z_1,z_2$ with $\alpha_m=\langle z_1,z_2\rangle$, 
and the rest of the proof follows from the similar argument of that in \cite[Corollary 2.9]{S}.
\end{proof}
\section{Real mutually unbiased bases}\label{MUB}
\begin{definition}
Let $M=\{M_i\}_{i=1}^f$ be a collection of orthonormal bases of $\mathbb{R}^d$.
$M$ is called real mutually unbiased bases (MUB) if any two vectors $x$ and $y$ from different bases
satisfy $\langle x,y \rangle =\pm 1/\sqrt{d}$.
\end{definition}
 
Let $M=\{M_i\}_{i=1}^f$ be a MUB, and put  $X=M \cup (-M)$. 
The angle set of $X$ is 
$$ A'(X)=\{\frac{1}{\sqrt{d}},0,-\frac{1}{\sqrt{d}},-1 \} .$$
We set 
$$\alpha_0=1,\quad\alpha_1=\frac{1}{\sqrt{d}},\quad\alpha_2=0,\quad \alpha_3=-\frac{1}{\sqrt{d}},\quad\alpha_4=-1 ,$$
and we define $R_k=\{(x,y)\in X\times X \mid \langle x,y\rangle=\alpha_k\}$.
Then $(X,\{R_k\}_{k=0}^4)$ is a $Q$-polynomial association scheme which is both $Q$-antipodal and $Q$-bipartite in \cite[Theorem 4.1]{LMO}.

Conversely let $(X,\{R_k\}_{k=0}^4)$ be a $Q$-polynomial association scheme which is both $Q$-antipodal and $Q$-bipartite,  
then the image of the embedding into first eigenspace by primitive idempotent $E_1$ is  $M\cup (-M)$, 
where $M$ is mutually unbiased bases in \cite[Theorem 4.2]{LMO}.

Applying \cite[Theorem 4.8]{BI} to the above scheme for $i=j=1$ using the paramerters in \cite[Appendix]{LMO}, 
we obtain the inequality $f\leq \frac{d}{2}+1$.   
We call $M$ a maximal MUB if this upper bound is attained.

\begin{lemma}\label{LSDofMUB}
$(X,\{R_i\}_{i=0}^4) $ is a $Q$-polynomial association scheme which is both $Q$-antipodal and $Q$-bipartite with $f$ $Q$-antipodal classes of size $2d$.
Assume $f\geq 3$.
Then for $z\in X$ and $j=1,3$ $(R_j(z),\{R_i\cap(R_j(z)\times R_j(z))\}_{i=0}^3\})$ is a $Q$-polynomial association scheme which is $Q$-antipodal with 
$(f-1)$ $Q$-antipodal classes of size $d$ and $a_1^*=\frac{d}{f-1}-2$. 
\end{lemma}
\begin{proof}  
It was shown in \cite[Secton 5]{S} that $(X,\{R_i\}_{i=0}^4) $ is triply regular, in particular $R_j(z)$ carries 
an association scheme for any $z\in X$, $j\in\{1,3\}$.
Let $X_j(z)$ be a derived design in $S^{d-2}$ of $X$ with respect to $z,\alpha_j$.
We verify the intersection numbers of $X_j(z)$. 
For $x,y\in X_j(z)$, we set 
$$p_{\alpha,\beta}(x,y)=|\{w\in X_j(z)\mid \langle x,w\rangle=\alpha,\langle w,y\rangle=\beta\}|.$$
The angle set of $X_j(z)$ is 
$$A(X_j(z))=\left\{\alpha_1:=\frac{\sqrt{d}-1}{d-1},\alpha_2:=\frac{-1}{d-1},\alpha_3:=\frac{-\sqrt{d}-1}{d-1}\right\},$$
$X_j(z)$ is a $s:=3$-distance set.
$X_j(z)$ is a $t:=2$-design in $S^{d-2}$, therefore $X_j(z)$ satisfies $t=2s-4$. 
And for any $\gamma=\langle x,y\rangle$, the intersection numbers $p_{\alpha_2,\alpha_2}(x,y)$, $p_{\alpha_2,\alpha_1}(x,y)$, $p_{\alpha_1,\alpha_2}(x,y)$ are 
independent of the choice of $x,y\in X_j(z)$ with $\gamma=\langle x,y\rangle$ as follows: 
\begin{align*}
p_{\alpha_2,\alpha_2}(x,y)=\begin{cases}
      0 &\text{ if }\ \langle x,y\rangle=\alpha_1, \\
      d-2& \text{ if }\ \langle x,y\rangle=\alpha_2, \\
      0 &\text{ if }\ \langle x,y\rangle=\alpha_3, 
\end{cases} p_{\alpha_2,\alpha_1}(x,y)=p_{\alpha_1,\alpha_2}(x,y)=\begin{cases}
     \frac{d+\sqrt{d}}{2}-1 & \text{ if }\ \langle x,y\rangle=\alpha_1, \\
     0 & \text{ if }\ \langle x,y\rangle=\alpha_2, \\
     \frac{d+\sqrt{d}}{2}   & \text{ if }\ \langle x,y\rangle=\alpha_3. 
\end{cases}
\end{align*}
For $0\leq\lambda\leq2$, $0\leq\mu\leq2$ and $(\lambda,\mu)\neq (1,2),(2,1),(2,2)$, 
we obtain a system of $6$ linear equations 
\begin{align*} \sum\limits_{\substack{1 \leq l \leq 3\\ 1 \leq m \leq 3\\ (l,m)\neq (2,2),(2,1),(1,2)}}\!\!\!\!\!\!\!\!\!\!\!\!\!\!\!\!\!\alpha_l^\lambda \beta_m^\mu p_{\alpha_l,\alpha_m}(x,y)  &=
   |X_j(z)|F_{\lambda,\mu}(\langle x,y\rangle)-\langle x,y\rangle^\lambda-\langle x,y\rangle^\mu-\alpha_2^\lambda \alpha_2^\mu p_{\alpha_2,\alpha_2}^j(x,y)\notag \\[-8mm] 
&\quad-\alpha_2^\lambda \alpha_1^\mu p_{\alpha_2,\alpha_1}^j(x,y)-\alpha_1^\lambda \alpha_2^\mu p_{\alpha_1,\alpha_2}^j(x,y), 
\end{align*}
where $F_{\lambda,\mu}(t)$ is defined in \cite[Section 7]{DGS}. 
$\{p_{\alpha_i,\alpha_j}(x,y)\mid 1\leq i,j\leq3, (i,j)\neq (2,2),(2,1),(1,2)\}$ is uniquely determined by Theorem~\ref{coherent}. 
The intersection matrices $B_i$ and the second eigenmatrix $Q$ are as follows:
\begin{align*}
B_1&=\begin{pmatrix}
0&1&0&0\\
\frac{(f-2)(d+\sqrt{d})}{2}&\frac{(f-3)(d+3\sqrt{d})}{4}&\frac{d+2\sqrt{d}}{4}&\frac{d+\sqrt{d}}{4}\\
0&\frac{d+\sqrt{d}-2}{2}&0&\frac{d+\sqrt{d}}{2}\\
0&\frac{(f-3)(d-\sqrt{d})}{4}&\frac{(f-2)d}{4}&\frac{(f-3)(d+\sqrt{d})}{4}
\end{pmatrix},\\
B_2&=\begin{pmatrix}
0&0&1&0\\
0&\frac{d+\sqrt{d}-2}{2}&0&\frac{d+\sqrt{d}}{2}\\
d-1&0&d-2&0\\
0&\frac{d-\sqrt{d}}{2}&0&\frac{d-\sqrt{d}-2}{2}
\end{pmatrix},\\
B_3&=\begin{pmatrix}
0&0&0&1\\
0&\frac{(f-3)(d-\sqrt{d})}{4}&\frac{(f-2)d}{4}&\frac{(f-3)(d+\sqrt{d})}{4}\\
0&\frac{d-\sqrt{d}}{2}&0&\frac{d-\sqrt{d}-2}{2}\\
\frac{(f-2)(d-\sqrt{d})}{2}&\frac{(f-3)(d-\sqrt{d})}{4}&\frac{(f-2)(d-2\sqrt{d})}{4}&\frac{(f-3)(d-3\sqrt{d})}{4}
\end{pmatrix},\\
Q&=\begin{pmatrix}
1&d-1&(f-1)(d-1)&f-1\\
1&\sqrt{d}-1&-\sqrt{d}+1&-1\\
1&-1&-f+1&f-1\\
1&-\sqrt{d}-1&\sqrt{d}+1&-1
\end{pmatrix},
\end{align*} 
and hence the Krein matrix $B_1^*$ is given as follows:
\begin{align*}
B_1^*&=\begin{pmatrix}
0&1&0&0\\
d-1&\frac{d}{f-1}-2&\frac{d}{f-1}&0\\
0&\frac{(f-2)d}{f-1}&\frac{(f-2)d}{f-1}-2&d-1\\
0&0&1&0
\end{pmatrix}.
\end{align*} 
Therefore $X_1(z)$ is a $Q$-polynomial association scheme which is $Q$-antipodal.
\end{proof}

The following Theorem shows that maximal MUB carries a quadruply regular association scheme 
and the quadruple regularity of an association scheme obtained from MUB is characterized by its parameter.  
\begin{theorem}
$(X,\{R_i\}_{i=0}^4) $ is a $Q$-polynomial association scheme which is both $Q$-antipodal and $Q$-bipartite.
Then the following conditions are equivalent:
\begin{enumerate}
\item $(X,\{R_i\}_{i=0}^4)$ is quadruply regular, 
\item $f=\frac{d}{2}+1$.
\end{enumerate}
\end{theorem}
\begin{proof}
(1)$\Rightarrow$(2):Assume $(X,\{R_i\}_{i=0}^4) $ is quadruply regular.
Then $X_1(z)$ is triply regular for any $z\in X$.
By Lemma~\ref{LSDofMUB} and Theorem~\ref{triLSD}, $\frac{d}{f-1}-2=0$.
Therefore $f=\frac{d}{2}+1$ holds.

(2)$\Rightarrow$(1):
By \cite[Corollary 5.3]{S} it is sufficient to show that the assumption of Corollary~\ref{quadruple} is satisfied.
 
(i) When  $\langle z_1,z_2\rangle=\alpha_2$, 
$\{(i,j)\mid 1\leq i\leq \frac{s-1}{2},1\leq j\leq s-1, p_{i,j}^m\neq0\} \cup\{(i,j)\mid \frac{s-1}{2}\leq i\leq \frac{s+1}{2},1\leq j\leq\frac{s+1}{2}, p_{i,j}^m\neq0\}$ is $\{(1,1),(1,3),(2,2)\}$.
$X_{i,j}^2=X_{i,j}^2(z_1,z_2)$ is a $3$-design in $S^{d-3}$ for $(i,j)\in \{(1,1),(1,3),(2,2)\}$.
Indeed $X_{2,2}^2$ is a cross polytope in $S^{d-3}$.
$|X_{1,1}^2|=p_{1,1}^2=\frac{d^2}{4}$, $|X_{1,3}^2|=p_{1,3}^2=\frac{d^2}{4}$ where $p_{1,1}^2$ and $p_{1,3}^2$ are the intersectioin numbers of $X$ in \cite[6 Appendix]{ABS}.
And the angle sets $A(X_{1,1}^2)=A(X_{1,3}^2)=\{\frac{\sqrt{d}-2}{d-2},\frac{-2}{d-2},\frac{-\sqrt{d}-2}{d-2}\}$ hold, so 
Gegenbauer polynomial expansion of their annihilator polynomial $F(t):=\prod_{\alpha\in A'(X_{i,j}^2)}\frac{t-\alpha}{1-\alpha}$ is
$$F(t)=\frac{4}{d^2}Q_0(t)+\frac{2(d^2+6)(d-2)}{d^3(d-1)}Q_1(t)+\frac{(d-2)^3(d+3)}{d^3(d-1)}Q_2(t)+\frac{6(d-2)(d-3)}{d^2(d-1)},$$
therefore $X_{1,1}^2$ and $X_{1,3}^2$ are $3$-designs in $S^{d-3}$ by \cite[Theorem 6.5]{DGS}.   
We renumber as follows:
$$X_1=X_{2,2}^2,\quad X_2=X_{1,1}^2,\quad X_3=X_{1,3}^2.$$
We define $s_{i,j}=|A(X_i,X_j)|$.
Then the matrix $(s_{i,j})$ is 
$$\begin{pmatrix}
1&2&2\\
2&3&3\\
2&3&3
\end{pmatrix}.$$
If $s_{i,j}+s_{j,k}-2\leq 3$, that is, when one of the $i,j,k$ at least is equal to $1$,
then the assumption (1) of Theorem~\ref{coherent} holds.

If $s_{i,j}+s_{j,k}-3=3$, that is, when 
\begin{align}\label{eq:3}(i,j,k)\in \{(l,m,n)\mid 2\leq l,m,n\leq 3\} ,\end{align}
And $X_2\cup X_3$ carries a subconstituent association scheme $R_1(z_1)$ of $X$ 
whose parameters are independent of $z_1$ by Lemma~\ref{LSDofMUB},
therefore those for $(2,3,3)$ (respectively $(2,3,2)$, $(3,2,3)$) are determined by those for $(2,2,3)$ (respectively $(2,2,2)$, $(3,3,3)$).
The intersection numbers $\{p_{\alpha,\beta}^j \mid \alpha=\alpha_{i,j}^2 \text{ or } \beta=\alpha_{j,k}^2\}$ for $x\in X_i$, $y\in X_k$ and $(i,j,k)\in \{(2,2,2),(3,3,3),(2,2,3)\}$ are given in Table~\ref{tb:t1}.  
These numbers are independent of $z_1,z_2\in X$ with $\langle z_1,z_2\rangle=\alpha_2$.
Hence the assumption of (2) of Theorem~\ref{coherent} holds for $i,j,k$ $(i,j,k)$ in (\ref{eq:3}).  

(ii) When  $\langle z_1,z_2\rangle=\alpha_1$, 
$\{X_{i,j}^m(z_1,z_2)\mid 1\leq i\leq \frac{s-1}{2},1\leq j\leq s-1, p_{i,j}^m\neq0\} \cup\{X_{i,j}^m(z_1,z_2)\mid \frac{s-1}{2}\leq i\leq \frac{s+1}{2},1\leq j\leq\frac{s+1}{2}, p_{i,j}^m\neq0\}$ is 
$\{X_{1,1}^1,X_{1,2}^1,X_{1,3}^1,X_{2,1}^1\}$.
$X_{i,j}^1=X_{i,j}^1(z_1,z_2)$ is a $2$-design in $S^{d-3}$  for $(i,j)\in \{(1,1),(1,2),(1,3),(2,1)\}$.
Indeed $X_{1,2}^1$, $X_{2,1}^1$ are regular simplexes in $S^{d-3}$.
And $X_{1,1}^1$ and $X_{1,3}^1$ are subconstituents of $X_1(z_1)$ with respect to $z_2\in X_1(z_1)$.
$X_1(z_1)$ is a $Q$-polynomial association scheme by Theorem~\ref{LSDofMUB} with $a_1^*=0$, so Lemma 4.2 in \cite{S} implies that $X_{1,1}^1$ and $X_{1,3}^1$ are $2$-designs in $S^{d-3}$.   
We renumber as follows:
\begin{align*}
X_1&=X_{2,1}^1,\quad X_2=X_{1,2}^1,\quad X_3=X_{1,1}^1,\quad X_4=X_{1,3}^1.
\end{align*}
We define $s_{i,j}=|A(X_i,X_j)|$.
Then the matrix $(s_{i,j})$ is 
$$\begin{pmatrix}
1&2&2&2\\
2&1&2&2\\
2&2&3&3\\
2&2&3&3
\end{pmatrix}.$$
If $s_{i,j}+s_{j,k}-2\leq 2$, that is, when  
$$(i,j,k)\in \{(l,m,n)\mid 1\leq m\leq2,1\leq l,n\leq4 \text{ or } 3\leq m\leq4,1\leq l,n\leq2 \},$$
then the assumption (1) of Theorem~\ref{coherent} holds.

If $s_{i,j}+s_{j,k}-3=2$, that is, when 
\begin{align}\label{eq:1}(i,j,k)\in \{(l,m,n)\mid 1\leq l\leq2, 3\leq m,n\leq4 \text{ or } 3\leq l,m\leq4,1\leq n\leq2\},\end{align}
or if $s_{i,j}+s_{j,k}-4=2$, that is, when
\begin{align}\label{eq:2}
(i,j,k)\in \{(l,m,n)\mid 3\leq l,m,n\leq4\},
\end{align}
we do not show that the $(i,j,k)$ in (\ref{eq:1}) (respectively (\ref{eq:2})) satisfy the assumption (2) (respectively (3)) of Theorem~\ref{coherent}, 
directly verify that the intersection numbers on $X_j$ for $x\in X_i$, $y\in X_k$ are independent of $x,y$ and of $z_1,z_2$ by using the triple regularity of subconsitituents of $X$. 
$X_2\cup X_3\cup X_4$ carries a subconstituent association scheme $R_1(z_1)$ 
which is obtained from a system of linked symmetric designs with $a_1^*=0$, 
and $X_2,X_3,X_4$ are the subconstituents of $R_1(z_1)$ with respect to $z_2\in R_1(z_1)$.
$R_1(z_1)$ is triply regular, so $X_2\cup X_3\cup X_4$ carries a coherent configuration whose parameters are independent of $z_2$. 
The parameters of $X_2\cup X_3\cup X_4$ depends on those of $R_1(z_1)$ which is independent of $z_1$. 
Therefore the parameters of $X_2\cup X_3\cup X_4$ are independent of $z_1,z_2$ with $\langle z_1,z_2\rangle=\alpha_1$.
Interchanging $z_1$ with $z_2$ and using $X_4=-X_{3,1}^1$, 
we can show  $X_1\cup X_3\cup X_4$ carries a coherent configuration whose parameters are independent of $z_1,z_2$ with $\langle z_1,z_2\rangle=\alpha_1$.

(iii) The case $\langle z_1,z_2\rangle=\alpha_3$ is similar to the case $\langle z_1,z_2\rangle=\alpha_1$.

By Corollary~\ref{quadruple},  we obtain the desired result.
\end{proof}

\begin{table}[h]
\begin{center}
\caption{the values of $p_{\alpha,\beta}^j(x,y)$, where $x\in X_i(z)$, $y\in X_k(z)$}
\begin{tabular}{|@{}c@{}|c|c||@{}c@{}|c|c|}\hline
$(i,j,k)$&$(\alpha,\beta)$&$p_{\alpha,\beta}^j(x,y)$&$(i,j,k)$&$(\alpha,\beta)$&$p_{\alpha,\beta}^j(x,y)$ \\ \hline
 & $(\alpha_{i,j}^2,\alpha_{j,k}^2)$ & 
$\left\{\begin{array}{ll}
     0 & \langle x,y\rangle=\alpha_{i,k}^1 \\
     \frac{d}{2}-1 & \langle x,y\rangle=\alpha_{i,k}^2 \\
     0 & \langle x,y\rangle=\alpha_{i,k}^3
\end{array}\right.$ &&$(\alpha_{2,2}^2,\alpha_{2,3}^2)$ &$\left\{\begin{array}{ll}
     0 & \langle x,y\rangle=\alpha_{2,3}^1 \\
     \frac{d}{2}-1 & \langle x,y\rangle=\alpha_{2,3}^2 \\
     0 & \langle x,y\rangle=\alpha_{2,3}^3
\end{array}\right.$ \\ \cline{2-3} \cline{5-6}
$\begin{matrix}(2,2,2)\\(3,3,3)\end{matrix}$ & $\begin{matrix}(\alpha_{i,j}^2,\alpha_{j,k}^1)\\(\alpha_{i,j}^1,\alpha_{j,k}^2)\end{matrix}$ & 
$\left\{\begin{array}{ll}
     \frac{d+2\sqrt{d}}{4}-1 &  \langle x,y\rangle=\alpha_{i,k}^1 \\
     0 & \langle x,y\rangle=\alpha_{i,k}^2 \\
     \frac{d+2\sqrt{d}}{4} & \langle x,y\rangle=\alpha_{i,k}^3
\end{array}\right.$ & & $(\alpha_{2,2}^2,\alpha_{2,3}^1)$ & $\left\{\begin{array}{ll}
     \frac{d}{4}-1 & \langle x,y\rangle=\alpha_{2,3}^1 \\
     0 & \langle x,y\rangle=\alpha_{2,3}^2 \\
     \frac{d}{4} & \langle x,y\rangle=\alpha_{2,3}^3
\end{array}\right.$ \\ \cline{2-3} \cline{5-6}
 & $\begin{matrix}(\alpha_{i,j}^2,\alpha_{j,k}^3)\\(\alpha_{i,j}^3,\alpha_{j,k}^2)\end{matrix}$ & 
$\left\{\begin{array}{ll}
     \frac{d-2\sqrt{d}}{4} & \langle x,y\rangle=\alpha_{i,k}^1 \\
     0 & \langle x,y\rangle=\alpha_{i,k}^2 \\
     \frac{d-2\sqrt{d}}{4}-1 & \langle x,y\rangle=\alpha_{i,k}^3
\end{array}\right.$&$(2,2,3)$ &$(\alpha_{2,2}^2,\alpha_{2,3}^3)$ &$\left\{\begin{array}{ll}
     \frac{d}{4}-1 & \langle x,y\rangle=\alpha_{2,3}^1 \\
     0 & \langle x,y\rangle=\alpha_{2,3}^2 \\
     \frac{d}{4} & \langle x,y\rangle=\alpha_{2,3}^3
\end{array}\right.$   \\ \cline{1-3} \cline{5-6}
&&&&$(\alpha_{2,2}^1,\alpha_{2,3}^2)$ &$\left\{\begin{array}{ll}
     \frac{d+2\sqrt{d}}{4} & \langle x,y\rangle=\alpha_{2,3}^1 \\
     0 & \langle x,y\rangle=\alpha_{2,3}^2 \\
     \frac{d+2\sqrt{d}}{4} & \langle x,y\rangle=\alpha_{2,3}^3
\end{array}\right.$ \\ \cline{5-6}
&&&&$(\alpha_{2,2}^3,\alpha_{2,3}^2)$ &$\left\{\begin{array}{ll}
     \frac{d-2\sqrt{d}}{4} & \langle x,y\rangle=\alpha_{2,3}^1 \\
     0 & \langle x,y\rangle=\alpha_{2,3}^2 \\
     \frac{d-2\sqrt{d}}{4} & \langle x,y\rangle=\alpha_{2,3}^3
\end{array}\right.$\\ \hline
\end{tabular}\label{tb:t1}
\end{center}
\end{table}

\begin{remark}
Let $M$ be a maximal MUB and $X=M\cup (-M)$.
It was already shown in \cite[Theorem 5]{A]BS} that $\{x\in X\mid \langle x,z_1\rangle=\langle x,z_2\rangle=\frac{1}{\sqrt{d}}\}$ for $z_1,z_2 \in X$ such that $\langle z_1,z_2\rangle=0$ carries an association scheme.
\end{remark}
\begin{remark}
One might wonder if integralities of intersection numbers of  the  above quadruply regular association scheme imply new necessary condition for existing maximal mutually unbiased bases, but these integralities show $d=\frac{k^2}{16}$, which is already known.  
\end{remark}
\section{Quadruple regularity of linked systems of symmetric designs}
Finally we consider whether the linked of symmetric designs with $a_1^*=0$ could become quadruply regular or not.
We denote the collection of all $k$-subsets of $\Omega$ by $\binom{\Omega}{k}$.
Let $(X_i,X_j,I_{i,j})$ be an incidence structure satisfying 
$X_i\cap X_j=\emptyset$, $I_{j,i}^t=I_{i,j}$ for any distinct integers $i,j \in \{1,\dots,f\}$. 
We put $X=\bigcup_{i=1}^f X_i$,
$I=\bigcup_{i\neq j} I_{i,j}$.
Let $(X,I)$ be a linked system of symmetric designs with $1<k<v-1$.
By \cite[Theorem 1]{C}, $n=k-\lambda$ is a square number.
Since $k<v-1$, we have $n\neq1$.
Hence $n\geq4$ and we have $v\geq15$.
We define 
$$\alpha(S)=|R_1(x_1)\cap R_1(x_2)\cap R_1(x_3)\cap R_1(x_4)|,$$
for $S=\{x_1,x_2,x_3,x_4\}\in \binom{X_1}{4}$. 
 
Counting in two ways the numbers of these sets 
\begin{align*}
&\{(S,y)\in \binom{X_1}{4} \times\bigcup\nolimits_{i=2}^fX_i\mid(x,y)\in R_1 \text{ for any }x\in S\},\\ 
&\{(S,T)\in \binom{X_1}{4} \times \binom{\bigcup\nolimits_{i=2}^fX_i}{2}\mid (x,y)\in R_1\text{ for any }x\in S,y\in T\},
\end{align*} 
we have the following equalities:
\begin{align}
\sum\limits_{S\in \binom{X_1}{4}}\alpha(S)&=(f-1)v\binom{k}{4},\\
\sum\limits_{S\in \binom{X_1}{4}}\binom{\alpha(S)}{2}&=\frac{1}{2}(f-1)v\left((f-2)k\binom{\sigma}{4}+(v-1)\binom{\lambda}{4}+(f-2)(v-k)\binom{\tau}{4}\right).
\end{align}
Using Cauchy-Schwarts inequality, we obtain 
$$ \frac{(k-1)k^2(v-k)^2(v-k-1)(v-2)}{(v-2k)(v-1)^2v}(k(v-3)(v-k)(v-2k)+\sqrt{n}(v-1)(v^2-6kv+v+6k^2))\geq0.$$
If equality holds, then we have $k(v-3)(v-k)(v-2k)+\sqrt{n}(v-1)(v^2-6kv+v+6k^2)=0$. 
Multiplying $k(v-3)(v-k)(v-2k)-\sqrt{n}(v-1)(v^2-6kv+v+6k^2)$, dividing by $vk(k-1)(v-k-1)(v-k)\neq0$, 
we obtain
$$v(v+1)^2+4k^2(v+3)-4kv(v+3)=0.$$  
Solving this quadratic equatuion by $k$, we obtain
$$k=\frac{v(v+3)\pm\sqrt{v(v-1)(v+3)}}{2(v+3)}.$$ 
Then it is necessary condition that $v(v-1)(v+3)$ is a square number.
The elliptic curve $y^2=v(v-1)(v+3)$ has rank $0$ and only $6$ integral points $(v,y)=(-2,6)$, $(-1,4)$, $(3,36)$, $(0,0)$, $(1,0)$, $(-3,0)$. 
It contradicts $v\geq15$.
Therefore $(X_1,\bigcup\nolimits_{i=2}^fX_i)$ does not become $4$-design.
Hence a linked system of symmetric designs does not carry a quadruply regular association scheme and
maximal MUB do not carry a quintuply regular association scheme.
However, let $(X,\{R_i\}_{i=0}^4)$ be a maximal MUB with $d>4$, the numbers 
\begin{align}\label{number}
|R_{i_1}(x_1)\cap R_{i_2}(x_2)\cap R_{i_3}(x_3)\cap R_{i_4}(x_4)\cap R_{i_5}(x_5)|
\end{align}
for $x_k\in X_1$ such that $(x_k,x_l)\in R_2$, $i_m\in\{1,3\}$ and $1\leq k,l,m\leq5$ are uniquely determined as follows. 
Let $M=\{M_1,\ldots,M_{d/2+1}\}$ be a maximal MUB.
Consider the orthogonal transformation on $M$ given by $M_1$ to standard basis.
Then the elements of $\bigcup_{i=2}^{\frac{d}{2}+1}X_i$ have form $\frac{1}{\sqrt{d}}(\pm1,\ldots,\pm1)$ since $M_1$ and $M_i$ are mutually unbiased. 
The binary code $\mathcal{C}$ is defined corresponding to the elements of $\bigcup_{i=2}^{\frac{d}{2}+1}X_i$ as follows:
$c=(c_i)\in\mathcal{C}$ corresponds to $x=(x_i)\in\bigcup_{i=2}^{\frac{d}{2}+1}X_i$, then $c_i=0,1$ according to whether $x_i=\frac{1}{\sqrt{d}},-\frac{1}{\sqrt{d}}$, respectively.
$\mathcal{C}$ is said to be a Kerdock-like code in \cite{ABS}.  
The weight enumerator of $\mathcal{C}$ is $W_{\mathcal{C}}(x,y)=x^d+\frac{d(d-2)}{2}x^{\frac{d+\sqrt{d}}{2}}y^{\frac{d-\sqrt{d}}{2}}+2(d-1)x^{\frac{d}{2}}y^{\frac{d}{2}}+\frac{d(d-2)}{2}x^{\frac{d-\sqrt{d}}{2}}y^{\frac{d+\sqrt{d}}{2}}+y^d$.
Then
\begin{align*}
W_{\mathcal{C}^\perp}(x,y)=\frac{1}{d^2}W_{\mathcal{C}}(x+y,x-y)
=x^d+\frac{d(d-1)(d-2)(d-4)}{360}x^{d-6}y^6+\cdots.
\end{align*}
Hence $\mathcal{C}$ is an orthgonal array whose strength is $5$ if $d>4$.
This implies that (\ref{number}) are uniquly determined for $x_k\in X_1$ such that $(x_k,x_l)\in R_2$, $i_m\in\{1,3\}$ and $1\leq k,l,m\leq5$.
\section*{Acknowledgements}
The author would like to thank Professor Akihiro Munemasa 
for helpful discussions.
This work is supported by Grant-in-Aid for JSPS Fellows. 

\end{document}